\documentclass[runningheads,orivec]{llncs}
\usepackage[T1]{fontenc}
\usepackage{color}
\usepackage[dvipsnames]{xcolor}
\usepackage{graphicx}%
\usepackage[utf8x]{inputenc}
\usepackage[english]{babel}

\usepackage{enumitem}
\setlist[itemize]{leftmargin=*}
\usepackage[toc,page]{appendix}
\usepackage{mathtools}
\usepackage{amsmath}
\usepackage{amsfonts}
\usepackage{amssymb}
\usepackage{animate}
\usepackage{multirow}
\usepackage{pdflscape}
\usepackage{todonotes}
\usepackage{verbatim}
\usepackage{float}
\usepackage[linesnumbered,ruled,vlined]{algorithm2e}
\usepackage{algpseudocode}

\newcommand{\n}{{\tfrac{\mu_2}{\mu_1}}}
\newcommand{\m}{{\tfrac{L_2}{\mu_1}}}
\NewDocumentCommand{\pospart}{ O{\ensuremath{\mu}} }{\ensuremath{[#1]_{+}}}
\NewDocumentCommand{\negpart}{ O{\ensuremath{\mu}} }{\ensuremath{[#1]_{-}}}
\usepackage{enumitem}
\setlist[itemize]{align=parleft,left=0pt..1em}
\usepackage{makecell}

\usepackage{pifont}
\usepackage{bm}
\usepackage{ascii}

\usepackage{comment}

\DeclareMathOperator*{\argmin}{argmin}
\DeclareMathOperator*{\argmax}{argmax}

\DeclareMathOperator*{\minimize}{minimize}

\DeclareMathOperator*{\dom}{dom}

\DeclareMathOperator*{\range}{range}

\usepackage{xparse}
\NewDocumentCommand{\Ei}{ O{k} O{x} }{\ensuremath{E_{#1} (#2)}}
\NewDocumentCommand{\TigLmu}{ O{i} O{\ensuremath{\gamma}} }{\ensuremath{T^{#1} (#2 L, #2 \mu)}}
\NewDocumentCommand{\Tinm}{ O{k} O{x} O{y} }{\ensuremath{T_{#1}(#2,#3)}}
\NewDocumentCommand{\alphaij}{ O{i} O{j} }{\ensuremath{\alpha_{[#1,#2]}}}

\usepackage{tikz,forest}
  \usetikzlibrary{shapes,arrows,fit,calc,positioning} %
  \tikzset{box/.style={draw, rectangle, thick, text centered, minimum height=0.5cm, minimum width=1cm}}
  \tikzset{line/.style={draw, thick, -latex'}}

\usepackage[unicode=true,psdextra,colorlinks]{hyperref}
\hypersetup{%
      bookmarksnumbered,
      breaklinks,
      hypertexnames = true,
      plainpages    = false,
      colorlinks = true,
      urlcolor   = blue,
      linkcolor  = MidnightBlue,
      filecolor  = Mulberry,      
      citecolor  = BrickRed,
      menubordercolor = BrickRed
    }

\usepackage[noabbrev, capitalise, nameinlink]{cleveref}
\newtheorem{assumption}{Assumption}
\Crefname{conjecture}{Conjecture}{Conjectures}
\Crefname{assumption}{Assumption}{Assumptions}
\Crefname{property}{Property}{Properties}
\Crefname{algorithm}{Algorithm}{Algorithms}

\begin{document}
\title{Tight Convergence Rates in Gradient Mapping for the Difference-of-Convex Algorithm}
\titlerunning{Tight Convergence Rates in Gradient Mapping for DCA}
%
\author{Teodor Rotaru\inst{1,2}\orcidID{0000-0003-2039-6228} \and
Panagiotis Patrinos\inst{1}\orcidID{0000-0003-4824-7697} \and
Fran\c{c}ois Glineur\inst{2}\orcidID{0000-0002-5890-1093}}
\authorrunning{Teodor Rotaru, Panagiotis Patrinos, Fran\c{c}ois Glineur}
%
\institute{Department of Electrical Engineering (ESAT-STADIUS), KU Leuven, Kasteelpark Arenberg 10, Leuven, 3001, Belgium
\email{ \{teodor.rotaru,panos.patrinos\}@esat.kuleuven.be} \and
Department of Mathematical Engineering (ICTEAM-INMA), UCLouvain, Av. Georges Lemaître 4, Louvain-la-Neuve, 1348, Belgium \\
\email{\{teodor.rotaru,francois.glineur\}@uclouvain.be}}
\maketitle              
\begin{abstract}
We establish new theoretical convergence guarantees for the difference-of-convex algorithm (DCA), where the second function is allowed to be weakly-convex, measuring progress via composite gradient mapping. Based on a tight analysis of two iterations of DCA, we identify six parameter regimes leading to sublinear convergence rates toward critical points and establish those rates by proving adapted descent lemmas.
We recover existing rates for the standard difference-of-convex decompositions of nonconvex-nonconcave functions, while for all other curvature settings our results are new, complementing recently obtained rates on the gradient residual. 
Three of our sublinear rates are tight for any number of DCA iterations, while for the other three regimes we conjecture exact rates, using insights from the tight analysis of gradient descent and numerical validation using the performance estimation methodology.
Finally, we show how the equivalence between proximal gradient descent (PGD) and DCA allows the derivation of exact PGD rates for any constant stepsize.
\end{abstract}
%
%
%

\section{Introduction}
We investigate the difference-of-convex (DC) formulation
$$
    \minimize_{x \in \mathbb{R}^d} F(x)  {}\coloneqq{}  f_1(x) - f_2(x),
$$
where $f_i \in \mathcal{F}_{\mu_i,L_i}$, $i=\{1,2\}$, denoting proper, lower semicontinuous (l.s.c.) functions $f_i:\mathbb{R}^d \rightarrow \mathbb{R}$ having curvature in the interval $[\mu_i, L_i]$, i.e., $\tfrac{L_i\|\cdot\|^2}{2}-f_i$ and $f_i - \frac{\mu_i \|\cdot\|^2}{2}$ are convex. If $f_i \in \mathcal{C}^2$, then $\mu_i$ and $L_i$ correspond to the minimum and maximum eigenvalues of the Hessian. When $L_i < \infty$, then $f_i$ is smooth. Moreover, $f_i$ is $\mu_i$-convex, including weakly-convexity under $\mu_i < 0$. We assume $\mu_1 \geq 0$ and allow $\mu_2 < 0$, which extends the standard DC case where both functions are convex. The difference-of-convex algorithm (DCA), stated in \cref{alg:dca}, is a standard method to solve the DC problem. Various work provide extensive analysis of DCA \cite{Dinh_Thi_dca_1997,DCA_Trust_region_1998,Horst_Thoai_DC_overview_1999,LeThi_2018_30_years_dev}, which is also referred as the convex-concave procedure (CCCP) \cite{CCCP_2001_init_Alan_Anand,CCCP_2009_Convergence,CCCP_Lipp_Boyd_2016}. In this work we focus on improving the theoretical guarantees. For the interested reader, a comprehensive list of DCA applications is found in the review paper \cite{LeThi_2018_30_years_dev}. 
\begin{algorithm}[t]
\caption{Difference-of-convex algorithm (DCA)}\label{alg:dca}
\KwData{$f_1 \in \mathcal{F}_{\mu_1,L_1}$, $f_2 \in \mathcal{F}_{\mu_2,L_2}$, with $\mu_1 \geq 0$, $\mu_1 + \mu_2 > 0$ or $\mu_1 = \mu_2 = 0$; $N \geq 1$ iterations starting from $x^0 \in \dom \partial f_1$}
\For{$ k = 0, \dots, N$}{
Select $g_2^k \in \partial f_2(x^k)$ \;
Select $x^{k+1} \in \argmin_{w\in \mathbb{R}^d} \{f_1(w) - \langle g_2^k, w \rangle\}$
}
\KwResult{Best iterate $x^k$ with $k = \argmin_{0 \leq k \leq N} \{\|x^{k} - x^{k+1}\|^2$\} }
\end{algorithm}%

If the entire sequence of iterates $\{x^k\}$ is bounded, then each accumulation point of it is a critical point of $F$ \cite[Theorem 3(iv)]{Dinh_Thi_dca_1997}. Inspired by this result, we aim to find the best iterate $x^k$ providing the lowest value of $\{\|x^{k} - x^{k+1}\|\}$ over $N$ iterations. This metric is referred as the \textit{best (composite) gradient mapping} \cite{Nesterov2013_Gradient_mapping_composite_functions,Nesterov_cvx_2018} or the \textit{prox-gradient mapping} \cite{Davis_Drusvyatskiy_2019_weakly_cvx}, and is complementary to the \textit{best residual gradient norm} $\{\|g_1^{k} - g_2^{k}\|\}$ over $N$ iterations, where $g_i \in \partial f_i(x^k)$, $i=\{1,2\}$, are any subgradients of $f_i$ evaluated at $x^k$. On this metric, \cite{abbaszadehpeivasti2021_DCA} provides the first tight rates for any convex functions $f_i$ when $F$ is nonconvex and nonconcave, using  the performance estimation problem (PEP) methodology \cite{drori_performance_2014,taylor_smooth_2017,Taylor_2017_SIAM_Composite_convex}. The analysis is later refined and extended in  \cite{rotaru2025tightanalysisdifferenceofconvexalgorithm} to accommodate weakly-convex functions $f_2$. For the standard DC setup, \cite{bregman_DCA_2023} provides convergence rates in terms of Bregman distance, viewing DCA as a instance of the Bregman proximal point algorithm, while  \cite{LeThi2018,yao2023_PL_DCA} provides linear rates under the Polyak-Łojasiewicz (PL) inequality. The Frank-Wolfe algorithm \cite{Yurtsever_Suvrit_FW_CCCP_2022} and the proximal gradient descent (PGD) \cite{LeThi_2018_30_years_dev} are also equivalent to DCA.%

\noindent\textit{Contributions.}
In \cref{thm:dca_rates_N_steps} we establish new theoretical convergence guarantees for DCA and PGD, measuring progress via composite gradient mapping. When both $f_1$ and $f_2$ are convex and $F$ is nonconvex-nonconcave, we refine the analysis from \cite{abbaszadehpeivasti2021_DCA}, providing standard descent lemmas in \cref{subsec:descent_lemmas}. Our results extend existing work by introducing new rates for all other curvature settings, complementing the best gradient residual rates from \cite{abbaszadehpeivasti2021_DCA,rotaru2025tightanalysisdifferenceofconvexalgorithm}. Building on insights from gradient descent (GD), we conjecture in \cref{conjecture:Regimes_p5_above_thr} the exact behavior of DCA across all curvature choices, with some regimes rigorously proven (\cref{thm:regime_p4_N_geq_3,thm:regime_p5_N_geq_3_dist_1}) and others remaining open. As shown in \cref{sec:PGD_rates}, the equivalence between PGD and DCA enables us to formulate precise PGD rates for any constant stepsize, including the ones exceeding the inverse Lipschitz constant. Finally, we propose a general adaptive curvature shifting procedure for DCA, leveraging PGD stepsize schedules aligned with established GD strategies. 

Our derivations are aided by the PEP technique, allowing to establish the tightness of the regimes and validate our conjecture. The setup of PEP for DCA is described in \cite{abbaszadehpeivasti2021_DCA,rotaru2025tightanalysisdifferenceofconvexalgorithm} and can be solved using dedicated software packages \cite{PESTO,PESTO_python}.

\smallskip
\noindent\textit{Theoretical setup and running assumptions.} 
Let $f$ be a proper, convex, l.s.c. function. Its subdifferential at $x \in \mathbb{R}^d$ is 
$\partial f(x) \coloneqq \{ g \in \mathbb{R}^d \,|\, f(y) {}\geq{} f(x) + \langle g, y-x \rangle, \forall y\in \mathbb{R}^d\}$. If $f$ is weakly-convex with $\mu<0$, let $\tilde{f}(x) \coloneqq f(x) - \mu \tfrac{\|x\|^2}{2}$ be a convex function. Since its subdifferential $\partial \tilde{f}(x)$ is well-defined, we take $\partial f(x) \coloneqq \partial \tilde{f}(x) + \mu x$ \cite[Proposition 6.3]{Bauschke_generalized_monotone_operators_2021}, which is well-defined. For $f$ differentiable at $x$ we have $\partial f(x) = \{ \nabla f(x)\}$. %
The domain and range of $f$ are $\dom f {}\coloneqq{} \{ x \in \mathbb{R}^d: f(x) < \infty \}$ and $\range f {}\coloneqq{} \{y \in \mathbb{R}: \exists \text{ } x \in \dom f \text{ with } y = f(x)\}$. The domain and range of the subdifferential are $\dom \partial f {}={} \{ x \in \mathbb{R}^d: \partial f(x) \neq \emptyset\}$ and $\range \partial f {}={} \cup \{\partial f(x): x \in \dom \partial f \}$. The convex conjugate of $f$ l.s.c. function is defined as $f^*(y) {}\coloneqq{} \sup_{x \in \dom f} \{ \langle y , x \rangle - f(x) \}$, where $f^*$ is closed and convex.%
\begin{assumption}\label{assumption:curvatures_and_well_def}
The objective function $F=f_1 - f_2$ is lower bounded by $F_{lo} \coloneqq \inf_x F$, where $f_1 \in \mathcal{F}_{\mu_1,L_1}$ and $f_2 \in \mathcal{F}_{\mu_2,L_2}$ with $0 
\leq \mu_1 < L_1 \leq \infty$, $\mu_2 \in (-\infty,\infty)$ and $L_2 \in (\mu_2, \infty]$. Additionally, assume $\varnothing \neq \dom \partial f_1 \subseteq \dom \partial f_2$ and $\range \partial f_2\subseteq \range \partial f_1$.
\end{assumption}%
Note that $F \in \mathcal{F}_{\mu_1 - L_2 \,,\, L_1 - \mu_2}$. The second part of \cref{assumption:curvatures_and_well_def} implies that the DCA iterations are well-defined, meaning there exists a sequence $\{x^k\}$, starting from $x^0 \in \dom \partial f_1$, generated by \cref{alg:dca}, such that $x^{k+1} \in \partial f_1^*(\partial f_2(x^k))$. By definition, this condition implies the existence of $g_1^{k+1} \in \partial f_1(x^{k+1})$ such that $g_1^{k+1} = g_2^{k}$, where $g_2^{k} \in \partial f_2(x^{k})$. This is the only characterization of DCA iterations employed in our derivations.%
DCA finds critical points $x^*$ satisfying $\partial f_2(x^*) \cap \partial f_1(x^*) \neq \emptyset$. These are stationary only if $f_2$ is smooth \cite[Exercise 10.10]{RockWets98}. If only $f_1$ is smooth, then only the inclusion $\partial (-f_2) (x^*) \subseteq -\partial f_2(x^*)$ is guaranteed \cite[Corollary 9.21]{RockWets98}, hence $\partial F(x^*) \subseteq \nabla f_1(x^*) - \partial f_2(x^*)$ and critical points are not necessary stationary. %
%

We assume that at each iteration a subgradient from $\partial f_1^*$ is available and only focus on the progress of the iterations $\min_{0\leq k \leq N} \{\|x^k - x^{k+1}\|^2\}$ subject to the initial condition $F(x^0)-F(x^{N+1}) \leq \Delta$, with $\Delta > 0$. To remove the dependence of the initial condition on the last iterate, in all our rates the initial objective gap can be further upper bounded by $F(x^0) - F_{lo}$.%
\section{Sublinear convergence rates}\label{sec:sublinear_convergence_rates}
\begin{theorem}\label{thm:dca_rates_N_steps}
Let $f_1 \in \mathcal{F}_{\mu_1,L_1}$ and $f_2 \in \mathcal{F}_{\mu_2,L_2}$ satisfying \cref{assumption:curvatures_and_well_def} and $\mu_1 + \mu_2 > 0$. Then after $N+1 \geq 2$ iterations of DCA starting from $x^0$ it holds
\begin{align}\label{eq:rate_N_steps_no_F*}
    \tfrac{1}{2} \min_{0\leq k \leq N} \{\|x^k - x^{k+1}\|^2\} {}\leq{}
    \frac{F(x^0)-F(x^{N+1})}{(\mu_1 + \mu_2) +  p_i(\mu_1,L_1, \mu_2,L_2) N},
\end{align}
where the expressions of the regimes $p_{i=\{1,\dots,6\}}$ are given in \cref{tab:DCA_regimes_one_step}.
\end{theorem}%
\begin{table}[thb!]
\centering
\caption{Coefficients $p_i$ from \cref{thm:dca_rates_N_steps}. Their domains cover all parameter space under the condition $\mu_1 + \mu_2 > 0$. Regimes $p_{\{1,2,3\}}$ are tight for any number of iterations. (s.c.: strongly convex; n.c.: nonconvex; ccv.: concave; n.ccv.: nonconcave)}%
\label{tab:DCA_regimes_one_step}
\begin{tabular}{|l|c|c|c|}
\hline
\textbf{\makecell[c]{Regime}} & \multicolumn{2}{c|}{\textbf{Domain}}  & \textbf{Description} \\ [2pt] \hline
$p_1 = \mu_1 + \mu_2 + \frac{(\mu_1-\mu_2)^2}{L_2-\mu_2}$  & $\mu_1 \geq \mu_2 \geq 0$ & $L_2 > \mu_1$
& \makecell[c]{ $f_1$, $f_2$ convex \\ $F$ n.c.-n.ccv.} \\ [2pt] \hline
$p_2 = \mu_1 + \mu_2 + \frac{(\mu_1-\mu_2)^2}{L_1-\mu_1}$  & $\mu_2 \geq \mu_1 \geq 0$
& $L_1 > \mu_2$
& \makecell[c]{ $f_1$, $f_2$ convex \\ $F$ n.c.-n.ccv.}  \\ [2pt] \hline
$p_3 =  \frac{(\mu_1+\mu_2)(L_2 + \mu_1)}{L_2 + \mu_2}  $  & $\mu_2 \in [\frac{-L_2 \mu_1}{L_2+\mu_1} ,0)$ &  $L_2 > \mu_1 > 0$
& \makecell[c]{ $f_1$ s.c., $f_2$ n.c. \\ $F$ n.c.-n.ccv.} \\ [2pt] \hline
$p_4 = \frac{\mu_1^2 (L_2 + \mu_1)}{L_2^2} $  & $\mu_2 \in \big[\frac{-L_2 \mu_1}{L_2+\mu_1}, L_2 \big)$ & $L_2 \in [0, \mu_1]$
& \makecell[c]{ $f_1$ s.c., $f_2$ n.c. \\ $F$ convex} \\ [2pt] \hline
$p_5 = \frac{\mu_1^2 (\mu_1 + \mu_2)}{\mu_2^2} $  & \multicolumn{2}{c|}{ \makecell{$\mu_2 < 0$, $\mu_2 \in \big({-\mu_1}, \frac{-L_2 \mu_1}{L_2+\mu_1}\big]$}}
& \makecell[c]{ $f_1$ s.c., $f_2$ n.c. \\ $F$ n.ccv. } \\ [2pt] \hline
$p_6 = \frac{\mu_2^2 (L_1 + \mu_2)}{L_1^2} $  & $\mu_2 \geq \mu_1 \geq 0$ &  $L_1 \in (0, \mu_2]$
& \makecell[c]{ $f_1$ convex, $f_2$ s.c. \\ $F$ ccv. } \\ [2pt] \hline
\end{tabular}%
\end{table}%
\begin{figure}[!htb]
    \centering
    \includegraphics[width=.65\textwidth]{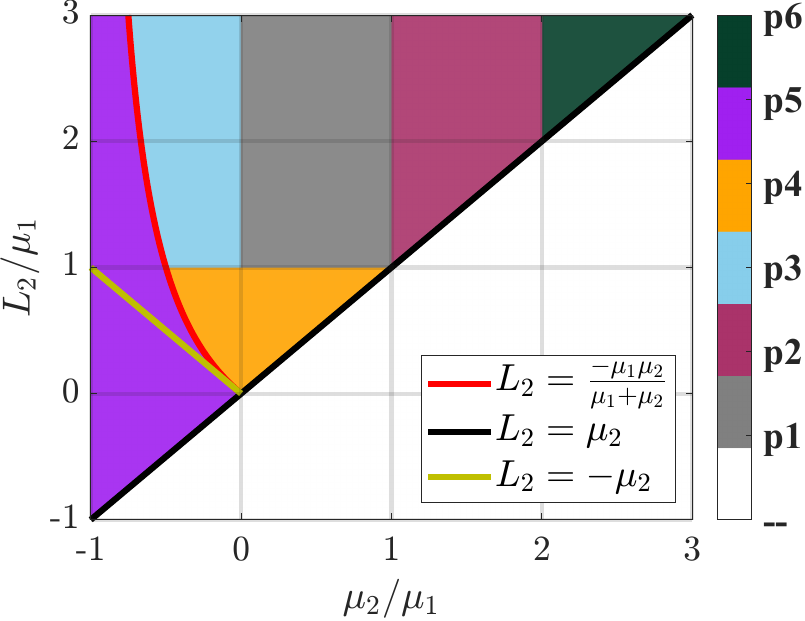}
    \caption{Domains of all tight regimes after two DCA iterations as proved in \cref{thm:dca_rates_N_steps}, given $\frac{L_1}{\mu_1} = 2$ and $\mu_1 > 0$. Regimes $p_1$, $p_2$ and $p_3$ are tight for any number of iterations $N \geq 2$. For the domain of $p_4$ we proved a linear rate for $N\geq 3$ (\cref{thm:regime_p4_N_geq_3}); same for the subdomain of $p_5$ satisfying the condition $L_2 + \mu_2 \leq 0$ (\cref{thm:regime_p5_N_geq_3_dist_1}). For $p_6$ and the rest of domain of $p_5$ we conjecture rates in \cref{subsec:conjecture_any_N}.
    }
    \label{fig:all_regimes_smooth_grd_map}
\end{figure}%

The proof is deferred to \cref{proof:thm:dca_rates_N_steps}. All regimes are exact for $N=\{0,1\}$ and their domains span the entire range of parameters satisfying the condition $\mu_1 + \mu_2 > 0$, which is sufficient for decreasing the objective (see \cref{prop:sufficient_decrease}). 

Based on extensive numerical experiments using PEP, we determine that regimes $p_{\{1,2,3\}}$ are tight for any number of iterations. Regimes $p_1$ and $p_2$ are derived in \cite[Proposition 3.1]{abbaszadehpeivasti2021_DCA}, by applying  Toland duality \cite{Toland1979} on the rates measuring the best gradient residual criterion. \cref{thm:dca_rates_N_steps} gives more explicit expressions and includes in the proofs a simpler derivation based on \cref{lemma:regime_p1,lemma:regime_p2}. These regimes assume $F$ is nonconvex-nonconcave ($L_1>\mu_2$ and $L_2>\mu_1$), with the active one indexed by $\argmax_i\{\mu_1,\mu_2\}$. All other regimes are new. \cref{fig:all_regimes_smooth_grd_map} shows the domains of all possible six regimes holding after two iterations. 

When $F$ is nonconvex-nonconcave, regime $p_3$ is active when $f_2$ is weakly-convex with $\mu_2 \in [\tfrac{-L_2 \mu_1}{L_2 + \mu_1},0)$. When $\mu_2$ is even lower, then regime $p_5$ is active. This is however not tight for $N \geq 2$, where a more complicated analysis is required. The exact rates expressions in this case are given in \cref{conjecture:Regimes_p5_above_thr}.

When $F$ is (strongly) convex ($\mu_1 \geq L_2$), the sublinear rates of regimes $p_4$ and $p_5$ are tight only for $N = 1$. These regimes are separated in this case by the threshold $\mu_2 = \frac{-L_2 \mu_1}{L_2 + \mu_1}$. When $N \geq 2$, we prove a linear rate for regime $p_4$ in \cref{thm:regime_p4_N_geq_3}, which is confirmed to be exact using PEP. Within the particular range of $p_5$ satisfying the condition $L_2 + \mu_2 \leq 0$, we demonstrate a linear rate in \cref{thm:regime_p5_N_geq_3_dist_1}. For the rest domain of $p_5$, we conjecture the tight rates in \cref{conjecture:Regimes_p5_above_thr}. Finally, regime $p_6$ characterizes the iterates over \textit{concave} objectives $F$ and is included for a complete covering of the parameter space.%
%
\subsection{Lemmas to prove sublinear rates} \label{subsec:descent_lemmas}
For compactness of the proofs, at each iteration $k$ we denote $\Delta x^{k} \coloneqq x^{k} - x^{k+1}$, $f_1^k=f_1(x^k)$, $f_2^k=f_2(x^k)$, $\Delta F(x^k) \coloneqq F(x^k) - F(x^{k+1})$, $g_1^k \in \partial f_1(x^k)$, $g_2^k \in \partial f_2(x^k)$, $G^k \coloneqq g_1^k - g_2^k$. Recall that $g_1^{k+1} = g_2^{k}$ for all DCA iterations.%
\begin{lemma}\label{lemma:nec_cycles_ineqs_dist1}
Let $f_1 \in \mathcal{F}_{\mu_1,L_1}$ and $f_2 \in \mathcal{F}_{\mu_2,L_2}$, with $L_2>\mu_2$ and $L_1>\mu_1 \geq 0$. Consider one DCA iteration connecting $x^k$ and $x^{k+1}$. Let $g_1^j \in \partial f_1(x^j)$, $g_2^j \in \partial f_2(x^j)$ and $G^j = g_1^j - g_2^j$, where $j=\{k,k+1\}$. The following inequalities hold:
    \begin{align}\label{eq:base_DeltaFk_dist_1}\tag{$B^{k}$}
    &\begin{aligned}        
    \Delta F(x^k)
        \geq &
    \tfrac{\mu_1 + \mu_2}{2} \|\Delta x^{k}\|^{2} {}+{}    
    \tfrac{\|G^{k+1} - \mu_2 \Delta x^{k}\|^2}{2(L_2-\mu_2)}  {}+{} 
    \tfrac{\|G^{k} - \mu_1 \Delta x^{k}\|^2}{2(L_1-\mu_1)} ; 
    \end{aligned} \\   
    & \begin{aligned}
    \langle G^{k}, \Delta x^{k} \rangle
        \geq
    \mu_1 \|\Delta x^{k}\|^{2} {}+{}
    \tfrac{\|G^{k} - \mu_1 \Delta x^{k}\|^2}{L_1-\mu_1} ;
\end{aligned}\label{eq:f1_dist_1_sum}  \tag{$C_{f_1}^{k}$}  \\
    & \begin{aligned}
    \langle G^{k+1}, \Delta x^{k} \rangle     
        \geq
    \mu_2 \|\Delta x^{k}\|^{2} {}+{}
    \tfrac{\|G^{k+1} - \mu_2 \Delta x^{k}\|^2}{L_2-\mu_2}.
\end{aligned}\label{eq:f2_dist_1_sum} \tag{$C_{f_2}^{k}$}
\end{align}%
\end{lemma}
\begin{proof}
Given $\mu < L$ and $f \in \mathcal{F}_{\mu,L}$, from \cite[Theorem 4]{taylor_smooth_2017}, \cite[Theorem 3.10]{Taylor_2017_SIAM_Composite_convex}, \cite[Theorem 3.2]{rotaru2024exactworstcaseconvergencerates} it holds $\forall x,y \in \mathbb{R}^d$, with $g^x \in \partial f(x)$, $g^y \in \partial f(y)$, that
\begin{align}\label{eq:Interp_hypoconvex}
    f(x)-f(y)-\langle g^y, x-y \rangle 
    \geq
    \tfrac{\mu}{2} \|x-y\|^{2} {}+{}
    \tfrac{1}{2(L-\mu)} \|g^x-g^y - \mu (x-y)\|^2. 
\end{align}%
We exploit the property $g_1^{k+1} = g_2^k$, implying 
$G^k = g_1^k-g_2^k = g_1^k - g_1^{k+1}$ and 
$G^{k+1} = g_1^{k+1}-g_2^{k+1} = g_2^k - g_2^{k+1}$.
By writing inequality \eqref{eq:Interp_hypoconvex} for: 
$f_1$ with the iterates $(x^{k},x^{k+1})$ and $(x^{k+1},x^{k})$ we get \eqref{eq:f_1_evals_in_xk_xk+1} and \eqref{eq:f_1_evals_in_xk+1_xk}, respectively;
$f_2$ with the iterates $(x^{k},x^{k+1})$ and $(x^{k+1},x^{k})$ we obtain \eqref{eq:f_2_evals_in_xk_xk+1} and \eqref{eq:f_2_evals_in_xk+1_xk}, respectively:%
\begin{align}\label{eq:f_1_evals_in_xk_xk+1}
    f_1^{k}-f_1^{k+1} - \langle g_1^{k+1}, \Delta x^{k} \rangle 
    {}\geq{} &
    \tfrac{\mu_1}{2} \|\Delta x^{k}\|^{2} {}+{}
    \tfrac{\|G^{k} - \mu_1 \Delta x^{k}\|^2}{2(L_1-\mu_1)}; \\
\label{eq:f_1_evals_in_xk+1_xk}
    f_1^{k+1}-f_1^{k} + \langle g_1^{k}, \Delta x^{k} \rangle 
    {}\geq{} &
    \tfrac{\mu_1}{2} \|\Delta x^{k}\|^{2} {}+{}
    \tfrac{\|G^{k} - \mu_1 \Delta x^{k}\|^2}{2(L_1-\mu_1)}; \\
\label{eq:f_2_evals_in_xk_xk+1}
    f_2^{k}-f_2^{k+1} - \langle g_2^{k+1}, \Delta x^{k}\rangle 
    {}\geq{} &
    \tfrac{\mu_2}{2} \|\Delta x^{k}\|^{2} {}+{}
    \tfrac{\|G^{k+1} - \mu_2 \Delta x^{k}\|^2}{2(L_2-\mu_2)};  \\
\label{eq:f_2_evals_in_xk+1_xk}
    f_2^{k+1}-f_2^k {}+{} \langle g_2^{k}, \Delta x^{k} \rangle 
    {}\geq{} &
    \tfrac{\mu_2}{2} \|\Delta x^{k}\|^{2} {}+{}
    \tfrac{\|G^{k+1} - \mu_2 \Delta x^{k}\|^2}{2(L_2-\mu_2)}.
\end{align}%
Then summing \eqref{eq:f_1_evals_in_xk_xk+1} and \eqref{eq:f_2_evals_in_xk+1_xk} gives \eqref{eq:base_DeltaFk_dist_1}; 
inequality \eqref{eq:f1_dist_1_sum} results by summing \eqref{eq:f_1_evals_in_xk_xk+1} and \eqref{eq:f_1_evals_in_xk+1_xk};
inequality \eqref{eq:f2_dist_1_sum} is obtained by summing \eqref{eq:f_2_evals_in_xk_xk+1} and \eqref{eq:f_2_evals_in_xk+1_xk}. \qed
\end{proof}

\begin{proposition}[Decrease after one iteration]\label{prop:sufficient_decrease}
   Let $f_1 \in \mathcal{F}_{\mu_1,L_1}$ and $f_2 \in \mathcal{F}_{\mu_2,L_2}$ satisfying \cref{assumption:curvatures_and_well_def}. Consider one DCA iteration connecting $x^{k}$ and $x^{k+1}$. Then it holds%
   \begin{align}\label{eq:standard_decrease}
   F(x^k) - F(x^{k+1})
        {}\geq{}
   \tfrac{\mu_1 + \mu_2}{2} \|x^k - x^{k+1}\|^2.
   \end{align}
   Consequently, the objective function $F$ decreases after each DCA iteration if $\mu_1 + \mu_2 \geq 0$ and the decrease is strict if $\mu_1 + \mu_2 > 0$, unless $x^{k+1}=x^{k}$.
\end{proposition}%
\begin{proof}
    Inequality \eqref{eq:standard_decrease} results by ignoring the mixed squares in \eqref{eq:base_DeltaFk_dist_1}. \qed
\end{proof}%
\cref{prop:sufficient_decrease} 
was first proved in \cite[Theorem 3, Proposition 2]{Dinh_Thi_dca_1997}.%
\begin{remark}[Proofs structure]\label{remark:proofs_structure}
    The proofs are based on multiplying the inequalities from \cref{lemma:nec_cycles_ineqs_dist1} with nonnegative weights. For regimes $p_{1,3,4,5}$ we sum:
    \begin{align}\label{eq:general_ineq_type}
        B^{k} {}+{} 
        \beta_1 \times C_{f_1}^{k+1} {}+{}
        \beta_2 \times C_{f_2}^{k},
    \end{align}
    while for the proofs of $p_{2,6}$ we sum:
    $
        B^{k+1} {}+{} 
        \beta_1 \times C_{f_1}^{k} {}+{}
        \beta_2 \times C_{f_2}^{k+1},
    $
    where $B^{\{k,k+1\}}$, $C_{f_1}^{\{k,k+1\}}$, $C_{f_2}^{\{k,k+1\}}$ are instances of inequalities \eqref{eq:base_DeltaFk_dist_1}, \eqref{eq:f1_dist_1_sum} and \eqref{eq:f2_dist_1_sum}, respectively. The weights $\beta_1$ and $\beta_2$ are functions of the curvature parameters.
\end{remark}
The demonstrations rewrite the inequalities from \cref{remark:proofs_structure} by building specific mixed squared norms of $G$ and $\Delta x$ at specific indices $k$, using \cref{property:identities_to_build_squares}. For space reasons, we only give the final inequalities containing all the squares.
\begin{property}\label{property:identities_to_build_squares}
    The following identities hold, where $i,j \in \{1,2\}$ and $G$ and $\Delta x$ are specialized when building the squares within the proofs: \\
    \smallskip
    $-\langle G, \Delta x \rangle {}={} \tfrac{1}{2\mu_{j}} \left[ \|G - \mu_{j} \Delta x\|^2 - \|G\|^2 - \|\mu_{j} \Delta x\|^2  \right]$ and \\
    \smallskip
    $\|G - \mu_{i} \Delta x\|^2 = \tfrac{\mu_{i}}{\mu_{j}} \|G - \mu_j \Delta x\|^2 - \tfrac{\mu_{i}}{\mu_{j}}(1-\tfrac{\mu_{i}}{\mu_{j}})  \|\mu_{j} \Delta x\|^2 {}+{} (1-\tfrac{\mu_{i}}{\mu_{j}}) \|G\|^2$.
\end{property}

\begin{lemma}[Regime $p_1$]\label{lemma:regime_p1}
   Let $f_1 \in \mathcal{F}_{\mu_1,L_1}$ and $f_2 \in \mathcal{F}_{\mu_2,L_2}$, with $L_2 > \mu_1 \geq \mu_2 \geq 0$. Then after two DCA iterations connecting $x^k$, $x^{k+1}$ and $x^{k+2}$ it holds
   $
        \Delta F(x^{k})
        {}\geq{}
   \big( \mu_1 + \mu_2 \tfrac{L_2 - \mu_1}{L_2 - \mu_2} \big) \tfrac{1}{2} \|\Delta x^{k}\|^2 + 
   \mu_1 \tfrac{\mu_1-\mu_2}{L_2-\mu_2}  \tfrac{1}{2} \|\Delta x^{k+1}\|^2.
   $
\end{lemma}
\begin{proof}
Let $\beta_2 = 0$ and $\beta_1 \coloneqq \frac{\mu_1 - \mu_2}{L_2 - \mu_2}$, which is nonnegative since $\mu_1 \geq \mu_2$. After replacing in \eqref{eq:general_ineq_type}, building squares and performing simplifications we get:%
    \begin{align*}
    \begin{aligned}
        \Delta F(x^{k})
      &  {}\geq{}
   \big( \mu_1 + \mu_2 \tfrac{L_2 - \mu_1}{L_2 - \mu_2} \big) \tfrac{1}{2} \|\Delta x^{k}\|^2 + 
   \mu_1 \tfrac{\mu_1-\mu_2}{L_2-\mu_2}  \tfrac{1}{2} \|\Delta x^{k+1}\|^2 + \\
   & \tfrac{\|G^{k} - \mu_1 \Delta x^{k}\|^2}{2(L_1-\mu_1)}  {}+{}
    \tfrac{\mu_2 \|G^{k+1} - \mu_1 \Delta x^{k}\|^2}{2\mu_1(L_2-\mu_2)}  {}+{} \tfrac{(L_1+\mu_1)(\mu_1-\mu_2) \|G^{k+1} - \mu_1 \Delta x^{k+1}\|^2}{2\mu_1(L_1-\mu_1)(L_2-\mu_2)}.
   \end{aligned}
   \end{align*}
Since $\mu_1 \geq \mu_2 \geq 0$, all coefficients of the squares are positive and the conclusion follows by dropping the mixed terms. \qed 
\end{proof}

\begin{lemma}[Regime $p_2$]\label{lemma:regime_p2}
   Let $f_1 \in \mathcal{F}_{\mu_1,L_1}$ and $f_2 \in \mathcal{F}_{\mu_2,L_2}$, with $L_1 > \mu_2 \geq \mu_1 \geq 0$. Then after two DCA iterations connecting $x^k$, $x^{k+1}$ and $x^{k+2}$ it holds
   $
   \Delta F(x^{k+1})
        {}\geq{}
    \mu_2 \tfrac{\mu_2-\mu_1}{L_1-\mu_1} \tfrac{1}{2} \|\Delta x^{k}\|^2 +
  \big( \mu_2 + \mu_1 \tfrac{L_1 - \mu_2}{L_1 - \mu_1} \big) \tfrac{1}{2} \|\Delta x^{k+1}\|^2.
   $
\end{lemma}
\begin{proof}
Let $\beta_1 = 0$ and $\beta_2 \coloneqq \frac{\mu_2 - \mu_1}{L_1 - \mu_1}$, nonnegative since $\mu_2 \geq \mu_1$, be the specific multipliers from \cref{remark:proofs_structure}. Building the mixed squares we obtain:%
$$
    \begin{aligned}
        & \Delta F(x^{k+1})
      {}\geq{}
   \big( \mu_2 + \mu_1 \tfrac{L_1 - \mu_2}{L_1 - \mu_1} \big) \tfrac{1}{2} \|\Delta x^{k+1}\|^2 + 
   \mu_2 \tfrac{\mu_2-\mu_1}{2(L_1-\mu_1)} \|\Delta x^{k}\|^2  {}+{}  \\
   & {}\qquad\quad{} \tfrac{\|G^{k+2} - \mu_2 \Delta x^{k+1}\|^2}{2(L_2-\mu_2)}  {}+{}
    \tfrac{\mu_1  \|G^{k+1} - \mu_2 \Delta x^{k+1}\|^2}{2\mu_2(L_1-\mu_1)} {}+{} 
   \tfrac{(L_2+\mu_2)(\mu_2-\mu_1) \|G^{k+1} - \mu_2 \Delta x^{k}\|^2}{2\mu_2(L_2-\mu_2)(L_1-\mu_1)}.
   \end{aligned}
$$%
Since $\mu_2 \geq \mu_1 \geq 0$, all coefficients of the squares are positive and the conclusion follows by dropping the mixed terms.\qed
\end{proof}%
%
%
\begin{property}[$T_1$]\label{property:T1}
    When $\mu_2 < 0$, we define $T_1 \coloneqq \mu_1 \big( \frac{L_2+\mu_1}{L_2^2} {}-{} \frac{\mu_1+\mu_2}{\mu_2^2} \big)$. One can check that $T_1 \geq 0$ if $\mu_2 \in (-\mu_1, \frac{-L_2 \mu_1}{L_2+\mu_1} \big]$ and $T_1 \leq 0$ if $\mu_2 \in [\frac{-L_2 \mu_1}{L_2+\mu_1}, L_2 \big)$.
\end{property}%
\begin{lemma}[Regime $p_3$]\label{lemma:regime_p3}
   Let $f_1 \in \mathcal{F}_{\mu_1,L_1}$ and $f_2 \in \mathcal{F}_{\mu_2,L_2}$, with $L_2 > \mu_1 > 0$, $\mu_2 \in [\frac{-L_2 \mu_1}{L_2+\mu_1}, 0 \big)$ and $\mu_1 + \mu_2 > 0$. Consider two DCA iterations connecting $x^k$, $x^{k+1}$ and $x^{k+2}$. Then
   $
   \Delta F(x^{k})
        {}\geq{}
   \big(\mu_1 + \tfrac{L_2 \mu_2}{L_2 + \mu_2}\big) \tfrac{1}{2}\|\Delta x^{k}\|^2 + \tfrac{\mu_1^2}{2(L_2+\mu_2)} \|\Delta x^{k+1}\|^2.
   $
\end{lemma}%
\begin{proof}
    Let $\beta_1 \coloneqq \frac{\mu_1}{L_2 + \mu_2}$ and $\beta_2 \coloneqq \frac{-\mu_2}{L_2 + \mu_2}$, which are nonnegative since $\mu_1 > -\mu_2 > 0$ and $L_2 + \mu_2 \geq \frac{-\mu_2 L_2}{\mu_1} > 0$. After replacing in \eqref{eq:general_ineq_type}, building squares and performing simplifications we obtain:%
    \begin{align*}
    \begin{aligned}
        \Delta F(x^{k})
        {}\geq{} &
   \big(\mu_1 + \tfrac{L_2 \mu_2}{L_2 + \mu_2}\big) \tfrac{1}{2} \|\Delta x^{k}\|^2 + \tfrac{\mu_1^2}{2(L_2+\mu_2)} \|\Delta x^{k+1}\|^2 {}+{} 
    \\ & \,
   \tfrac{ \|G^{k} - \mu_1 \Delta x^{k}\|^2}{2(L_1-\mu_1)} +
   \tfrac{(L_1 + \mu_1)  \|G^{k+1} - \mu_1 \Delta x^{k+1}\|^2}{2(L_1-\mu_1)(L_2+\mu_2)}.
    \end{aligned}
    \end{align*}
Since $L_2 > -\mu_2$, the coefficients of the mixed terms are positive and the conclusion follows by dropping them.
The coefficient of $\|\Delta x^k\|^2$ is positive if $T_1 \leq 0$ as it can be factorized as $\frac{L_2^2 \mu_2^2}{\mu_1(L_2^2 - \mu_2^2)} (-T_1)$. The inequality holds even for $T_1 > 0$.
    \qed
\end{proof}

\begin{lemma}[Regime $p_4$]\label{lemma:regime_p4}
   Let $f_1 \in \mathcal{F}_{\mu_1,L_1}$ and $f_2 \in \mathcal{F}_{\mu_2,L_2}$, with $L_2 \in [0, \mu_1]$, $\mu_2 \in \big[\frac{-L_2 \mu_1}{L_2+\mu_1}, L_2 \big)$ and $\mu_1 + \mu_2 > 0$. Consider two DCA iterations connecting $x^k$, $x^{k+1}$ and $x^{k+2}$. Then it holds
   $
   \Delta F(x^{k})
        {}\geq{}
   \tfrac{\mu_1^2(L_2 + \mu_1)}{L_2^2} \tfrac{1}{2} \|\Delta x^{k+1}\|^2.
   $
\end{lemma}
\begin{proof}
    Let $\beta_1 \coloneqq \frac{\mu_1(\mu_1+L_2)}{L_2^2}$ and $\beta_2 \coloneqq \frac{\mu_1}{L_2}$, both positive. After replacing in \eqref{eq:general_ineq_type}, building the mixed squares and performing simplifications we get:
    \begin{align*}
         \Delta F(x^{k})
        & {}\geq{}
   \tfrac{\mu_1^2(L_2 + \mu_1)}{L_2^2} \tfrac{1}{2} \|\Delta x^{k+1}\|^2 + \tfrac{\|G^k - \mu_1 \Delta x^{k}\|^2}{2(L_1-\mu_1)}  {}+{} 
   \\ 
   & \tfrac{(L_1 + \mu_1)(L_2 + \mu_1) \|G^{k+1} - \mu_1 \Delta x^{k+1}\|^2}{L_2^2(L_1 - \mu_1)}  +
   \tfrac{\mu_2^2 \, (-T_1) \|G^{k+1} - L_2 \Delta x^{k}\|^2}{\mu_1(L_2 - \mu_2)^2} .
    \end{align*}
    Since $T_1 \leq 0$ (\cref{property:T1}), the conclusion follows. \qed
\end{proof}

\begin{lemma}[Regime $p_5$]\label{lemma:regime_p5}
   Let $f_1 \in \mathcal{F}_{\mu_1,L_1}$ and $f_2 \in \mathcal{F}_{\mu_2,L_2}$, with $\mu_1 > -\mu_2 > 0$ such that $\mu_2 \in \big(-\mu_1, \frac{-L_2 \mu_1}{L_2+\mu_1}\big]$. Consider two DCA iterations connecting $x^k$, $x^{k+1}$ and $x^{k+2}$. Then it holds 
   $\Delta F(x^{k})
        {}\geq{}
   \tfrac{\mu_1^2(\mu_1 + \mu_2)}{\mu_2^2} \tfrac{1}{2} \|\Delta x^{k+1}\|^2$. %
\end{lemma}%
\begin{proof}
Let $\beta_1 \coloneqq \frac{\mu_1(\mu_1+\mu_2)}{\mu_2^2}$ and $\beta_2 \coloneqq \frac{\mu_1+\mu_2}{-\mu_2}$, both positive. After replacing in \eqref{eq:general_ineq_type}, building the mixed squares and performing simplifications we get:%
\begin{align*}   
\begin{aligned}
   \Delta F(x^{k})
       {}\geq{} &
   \tfrac{\mu_1^2(\mu_1 + \mu_2)}{2\mu_2^2} \|\Delta x^{k+1}\|^2
   + \tfrac{\|G^k - \mu_1 \Delta x^{k}\|^2}{2(L_1-\mu_1)} {}+{} 
   \\ 
   & \tfrac{(L_1 + \mu_1)(\mu_1 + \mu_2) \|G^{k+1} - \mu_1 \Delta x^{k+1}\|^2}{\mu_2^2(L_1 - \mu_1)}  +
   \tfrac{L_2^2 \, T_1 \|G^{k+1} - \mu_2 \Delta x^{k}\|^2}{\mu_1(L_2 - \mu_2)^2} .
\end{aligned}
\end{align*}%
Since $T_1 \geq 0$ (\cref{property:T1}), the conclusion follows. \qed
\end{proof}%
\begin{lemma}[Regime $p_6$]\label{lemma:regime_p6}
   Let $f_1 \in \mathcal{F}_{\mu_1,L_1}$ and $f_2 \in \mathcal{F}_{\mu_2,L_2}$, with $\mu_2 \geq \mu_1 \geq 0$ and $L_1 \in (0, \mu_2]$. Then after two DCA iterations connecting $x^k$, $x^{k+1}$ and $x^{k+2}$ we have $ \Delta F(x^{k+1})
        {}\geq{}
   \tfrac{\mu_2^2(L_1 + \mu_2)}{L_1^2} \tfrac{1}{2} \|\Delta x^{k}\|^2.
   $
\end{lemma}
\begin{proof}
    Let $\beta_2 \coloneqq \frac{\mu_2(\mu_2+L_1)}{L_1^2}$ and $\beta_1 \coloneqq \frac{\mu_2}{L_1}$, both positive. After replacing in \eqref{eq:general_ineq_type}, building the mixed squares and performing simplifications we get:
    \begin{align*}
    \Delta F(x^{k+1})
        {}\geq{} &
   \tfrac{\mu_2^2(L_1 + \mu_2)}{L_1^2} \tfrac{1}{2} \|\Delta x^{k}\|^2 + 
   \tfrac{\|G^{k+2} - \mu_2 \Delta x^{k+1}\|^2}{2(L_2-\mu_2)} {}+{} 
   \\ 
   & \, \tfrac{(L_1 + \mu_2)(L_2 + \mu_2) \|G^{k+1} - \mu_2 \Delta x^{k}\|^2}{L_1^2(L_2 - \mu_2)}  +
   \tfrac{\mu_1 \mu_2 + L_1(\mu_1+\mu_2) \|G^{k+1} - L_1 \Delta x^{k+1}\|^2}{L_1^2(L_1-\mu_1)} .
    \end{align*}
    The conclusion follows by dropping all mixed squares terms.\qed
\end{proof}
We remark the symmetry between regimes $p_1 \leftrightarrow p_2$ and $p_4 \leftrightarrow p_6$, as the latter ones reflect the behaviors of the former regimes, as the iterations would be applied in reverse order on the function $-F = f_2 - f_1$. In other words, \cref{lemma:regime_p2,lemma:regime_p6} are obtained from \cref{lemma:regime_p1,lemma:regime_p4}, respectively, by interchanging the curvatures $\mu_1 \leftrightarrow \mu_2$ and $L_1 \leftrightarrow L_2$ and the iterations order, i.e., from $k+2$ to $k$.%
\subsection{Proof of \cref{thm:dca_rates_N_steps}}\label{proof:thm:dca_rates_N_steps}
\begin{proof}[of \cref{thm:dca_rates_N_steps}]
    For regimes $p_1$, $p_3$, $p_4$ and $p_5$, by telescoping inequalities from \cref{lemma:regime_p1}, \cref{lemma:regime_p3}, \cref{lemma:regime_p4} and \cref{lemma:regime_p5}, respectively, for $k=\{0,\dots,N-1\}$ and taking the minimum over iterations differences norms:
    $$
        F(x^0)-F(x^{N}) {}\geq{} p_i(\mu_1,L_1, \mu_2,L_2) N \, \tfrac{1}{2}\min_{0 \leq k \leq N} \{\|x^k - x^{k+1}\|^2\}, \,\, i = \{1,3,4,5\},
    $$
    where $p_i$ denotes the sum of the two coefficients multiplying $\tfrac{1}{2}\|x^{k} - x^{k+1}\|^2$ and $\tfrac{1}{2}\|x^{k+1} - x^{k+2}\|^2$ in each inequality. Writing \eqref{eq:standard_decrease} for the iterates $x^{N}$ and $x^{N+1}$, we have $F(x^N) - F(x^{N+1}) {}\geq{} \frac{\mu_1+\mu_2}{2} \|x^{N} - x^{N+1}\|^2 $. Then the rate results by summing it to the previous inequality and taking the minimum in the right-hand side. For regimes $p_2$ and $p_6$, by telescoping inequalities from \cref{lemma:regime_p2,lemma:regime_p6}, respectively, for $k=\{1,\dots,N\}$ and taking the minimum over iterations:
    $$
        F(x^1)-F(x^{N+1}) {}\geq{} p_i(\mu_1,L_1, \mu_2,L_2) N \, \tfrac{1}{2}\min_{0 \leq k \leq N} \{\|x^k - x^{k+1}\|^2\}, \,\, i=\{2,6\}.
    $$
    Writing \eqref{eq:standard_decrease} for the iterates $x^{0}$ and $x^{1}$, summing it up and taking the minimum in the right-hand side leads to the rate. \qed
\end{proof}%
\section{Tight convergence rates for any number of iterations} \label{sec:tight_rates_any_N}%
Within this section we formulate the tight rates for any number of iterations for the regimes where the analysis for $N=2$, leading to sublinear rates, does not extend exactly for $N>2$, namely for regimes $p_{4,5,6}$. \\
\noindent \textbf{Notation.} We define $E_k: \mathbb{R}_{>0} \rightarrow \mathbb{R}_{>0}$, $\Ei[k][z]{}\coloneqq{}\sum_{j=1}^{2k} z^{-j}$, which is $\frac{-1+z^{-2k}}{1-z}$, if $z \neq 1$, and $2k$ is $z=1$. Note that $\Ei[0][z] = 0$. We denote $\pospart[z] \coloneqq \max\{0, z\}$. %
\subsection{Subdomains with proved rates}
\cref{thm:regime_p4_N_geq_3} provides the exact rate in the area of regime $p_4$, whereas \cref{thm:regime_p5_N_geq_3_dist_1} gives the exact behaviour in the subdomain of $p_5$, bounded by the condition $L_2 + \mu_2 \leq 0$. In both cases, the proofs only involve combining inequalities on consecutive iterations given in \cref{lemma:descent_lemma_stepsizes_geq_1}.%
\begin{theorem}\label{thm:regime_p4_N_geq_3}
   Let $f_1 \in \mathcal{F}_{\mu_1,L_1}$ and $f_2 \in \mathcal{F}_{\mu_2,L_2}$, with $L_2 \in [0, \mu_1]$, $\mu_1 + \mu_2 > 0$ and $\mu_2 \in \big[\frac{-L_2 \mu_1}{L_2+\mu_1}, L_2 \big)$. Then after $N+1 \geq 2$ iterations of DCA starting from $x^0$ it holds: 
   $
       \tfrac{1}{2} \|x^{N} - x^{N+1}\|^2 {}\leq{} \frac{F(x^0) - F(x^{N+1})}{\mu_1 + \mu_2 + \mu_1 \Ei[N][\frac{L_2}{\mu_1}]}.
   $
\end{theorem}%
\begin{proof}
We take $\beta_1 = \Ei[k+1][\m]$ and $\beta_2 = -1 + \m \Ei[k+1][\m]$ in \eqref{eq:general_ineq_type}, where $k=0,\dots,N-1$. After building the squares and simplifications we get:%
\begin{align*}
    \Delta F(x^{k})
        & \geq
    -\Ei[k][\m] \tfrac{\mu_1}{2} \|\Delta x^{k}\|^2 {}+{}
     \Ei[k+1][\m] \tfrac{\mu_1}{2} \|\Delta x^{k+1}\|^2 {}+{}
     \tfrac{\|G^k - \mu_1 \Delta x^{k}\|^2}{2(L_1-\mu_1)}  +
   \\ 
   & \,
   \tfrac{(L_1+\mu_1) \Ei[k+1][\frac{L_2}{\mu_1}] \|G^{k+1} - \mu_1 \Delta x^{k+1}\|^2 }{2\mu_1(L_1-\mu_1)}  +
   \tfrac{[ -1+ \frac{L_2+\mu_2}{\mu_1} \Ei[k+1][\frac{L_2}{\mu_1}] ] \, \|G^{k+1} - L_2 \Delta x^{k}\|^2}{2(L_2-\mu_2)} .
    \end{align*}%
We have $-1+ \frac{L_2+\mu_2}{\mu_1} \Ei[k+1][\m] {}\geq{} -1+ \frac{L_2+\mu_2}{\mu_1} \Ei[1][\m] {}={} \frac{\mu_2^2}{\mu_1(L_2-\mu_2)} (-T_1) {}\geq{} 0$, since $\Ei[k+1][\m]$ is  increasing with $k$, thus the mixed terms can be dropped. Telescoping for $k=0,\dots,N-1$ and adding \eqref{eq:base_DeltaFk_dist_1} with $k=N$ gives the conclusion. \qed
\end{proof}%
\begin{theorem}\label{thm:regime_p5_N_geq_3_dist_1}
   Let $f_1 \in \mathcal{F}_{\mu_1,L_1}$ and $f_2 \in \mathcal{F}_{\mu_2,L_2}$, with $\mu_1 + \mu_2 > 0$ and $L_2 + \mu_2 \leq 0$. Then after $N+1 \geq 2$ iterations of DCA starting from $x^0$ it holds %
   $
       \tfrac{1}{2} \|x^{N} - x^{N+1}\|^2 {}\leq{} \tfrac{F(x^0) - F(x^{N+1})}{\mu_1 + \mu_2 + \mu_1 \Ei[N][\frac{\mu_2}{\mu_1}]}.
   $
\end{theorem}%
\begin{proof}
We take $\beta_1 = \Ei[k+1][\n]$ and $\beta_2 = -\n \Ei[k+1][\n]$ in \eqref{eq:general_ineq_type}, where $k=0,\dots,N-1$. After building the squares and simplifications we get:
\begin{align*}
     \Delta F(x^{k})
        \geq &
    -\Ei[k][\n] \tfrac{\mu_1}{2} \|\Delta x^{k}\|^2 {}+{}
     \Ei[k+1][\n] \tfrac{\mu_1}{2} \|\Delta x^{k+1}\|^2+ 
     \tfrac{\|G^k - \mu_1 \Delta x^{k}\|^2}{2(L_1-\mu_1)} {}+{} 
   \\ 
   & \,
   \tfrac{(L_1+\mu_1) \Ei[k+1][\frac{\mu_2}{\mu_1}] \|G^{k+1} - \mu_1 \Delta x^{k+1}\|^2}{2\mu_1(L_1-\mu_1)}  {}+{} 
   \tfrac{[1 - \frac{L_2+\mu_2}{\mu_1} \Ei[k+1][\frac{\mu_2}{\mu_1}]] \|G^{k+1} - \mu_2 \Delta x^{k}\|^2}{2(L_2-\mu_2)}.
    \end{align*}
$\Ei[k+1][\n]$ is positive and monotonically increasing with index $k$ as $\mu_1 + \mu_2 > 0$. Then due to $L_2 + \mu_2 \leq 0$, we can drop all mixed squares. Telescoping for  $k=0,\dots,N-1$ and adding \eqref{eq:base_DeltaFk_dist_1} with $k=N$ gives the conclusion. \qed
\end{proof}%
\subsection{Conjectures for any number of iterations}\label{subsec:conjecture_any_N}%
For the domain of regime $p_5$ fulfilling the condition $L_2 + \mu_2 > 0$ we conjecture the exact behavior in \cref{conjecture:Regimes_p5_above_thr}.
\begin{conjecture}\label{conjecture:Regimes_p5_above_thr}
Let $f_1 \in \mathcal{F}_{\mu_1,L_1}$ and $f_2 \in \mathcal{F}_{\mu_2,L_2}$, with $\mu_1 \geq 0$ and $\mu_1 + \mu_2 > 0$. Consider $N+1 \geq 2$ iterations of DCA starting from $x^0$. \\
If $\mu_2 < 0 < \mu_1 < L_2$, thus $F$ is nonconvex-nonconcave, then:
\begin{align}\label{eq:regime_p5_p3_nonconvex_F}
    \tfrac{\mu_1}{2 } \min_{0\leq k \leq N} \left\{\|x^{k} - x^{k+1}\|^2\right\} {}\leq{}
    \tfrac{F(x^0)-F(x^{N+1})}{
    1 + \frac{\mu_2}{\mu_1} {}+{}
        \min \big\{ 
        P_N( \frac{L_2}{\mu_1}, \frac{\mu_2}{\mu_1} ) \,,\,
        \Ei[N][\frac{\mu_2}{\mu_1}]
        \big\}
    },
\end{align}
where
$
    P_N\left(\eta, \rho \right)
        {}\coloneqq{}
    \tfrac{(1 + \eta)(1+\rho)}{\eta+\rho} 
    \big(
        N + 
        \tfrac{(1-\eta)(1-\rho)}{\eta-\rho} 
        \sum_{k=0}^N 
        [
            \Ei[k][\eta] - \Ei[k][\rho]
        ]_{+}
    \big).
$ \\
If $L_2 \leq \mu_1$ and $\mu_1 > 0$, thus $F$ is (strongly) convex, then:
    \begin{align}\label{eq:regime_p5_p7_rate_strg_cvx}
    \tfrac{\mu_1}{2 } \min_{0\leq k \leq N} \left\{\|x^{k} - x^{k+1}\|^2\right\} {}\leq{}
    \tfrac{F(x^0)-F(x^{N+1})}{
    1 + \frac{\mu_2}{\mu_1} {}+{}
        \min \big\{ 
        \Ei[N][\frac{L_2}{\mu_1}] \,,\,
        \Ei[N][\frac{\mu_2}{\mu_1}]
        \big\}
    }.
\end{align}%
Finally, if $0 \leq \mu_1 < L_1 \leq \mu_2 < L_2$, thus $F$ is (strongly) concave, then:
    \begin{align}\label{eq:regime_p6_rate_strg_concave}
    \tfrac{\mu_2}{2 } \min_{0\leq k \leq N} \left\{\|x^{k} - x^{k+1}\|^2\right\}{}\leq{}
    \tfrac{F(x^0)-F(x^{N+1})}{
    1 + \frac{\mu_1}{\mu_2} {}+{}
        \min \big\{ 
        \Ei[N][\frac{L_1}{\mu_2}] \,,\,
        \Ei[N][\frac{\mu_1}{\mu_2}]
        \big\}
    }.
\end{align}%
\end{conjecture}%
For $N=1$, \cref{conjecture:Regimes_p5_above_thr} includes the rates corresponding to regimes $p_3$ and $p_5$ in \eqref{eq:regime_p5_p3_nonconvex_F} when $F$ is nonconvex, of regimes $p_4$ and $p_5$ in \eqref{eq:regime_p5_p7_rate_strg_cvx} when $F$ is (strongly) convex and of regime $p_6$ in \eqref{eq:regime_p6_rate_strg_concave} when $F$ is concave. The rate in \eqref{eq:regime_p5_p7_rate_strg_cvx} includes the ones proved in \cref{thm:regime_p4_N_geq_3} (within the domain of regime $p_4$) and in \cref{thm:regime_p5_N_geq_3_dist_1}. The expressions of the rates \eqref{eq:regime_p5_p3_nonconvex_F} and \eqref{eq:regime_p5_p7_rate_strg_cvx} are independent of $L_1$, hence they also hold when $f_1$ is \textbf{\textit{nonsmooth}}, i.e., $L_1 = \infty$. 
Note that the normalization of the iterates differences by $\mu_1$ or $\mu_2$ results in a homogeneous expression of the right-hand side denominators in the ratios $\frac{L_2}{\mu_1}$ and $\frac{\mu_2}{\mu_1}$, or $\frac{L_1}{\mu_2}$ and $\frac{\mu_1}{\mu_2}$, respectively.

\cref{conjecture:Regimes_p5_above_thr} is based on the exact rates developed for gradient descent (GD) derived in \cite{rotaru2024exactworstcaseconvergencerates} and is confirmed by PEP numerical experiments. GD is a particular case of the proximal gradient descent (PGD), which is iteration equivalent with DCA when $f_2$ is smooth.

\section{Deriving rates for proximal gradient descent (PGD)}\label{sec:PGD_rates}
Assume $F=\varphi + h$, with $\varphi \in \mathcal{F}_{\mu_{\varphi}, L_{\varphi}}$ smooth, $L_{\varphi} \in (0,\infty)$, $\mu_{\varphi} \in (-\infty, L_{\varphi})$, and $h \in \mathcal{F}_{\mu_{h}, L_{h}}$ proper, closed, l.s.c. and convex, such that $\mu_h \geq 0$ and $L_h \in (0,\infty]$. %

The PGD iteration with stepsize $\gamma > 0$, starting from $x^k$, reads $x^{k+1} = \argmin_{w \in \mathbb{R}^d} \{h(w) + \tfrac{1}{2\gamma} \|w-x^{k} + \gamma \nabla \varphi(x^{k})\|^2 \}$. This iteration is exactly the DCA one applied to $f_1 = h + \tfrac{\|\cdot\|^2}{2\gamma}$ and $f_2 = \tfrac{\|\cdot\|^2}{2\gamma} - \varphi$, as shown in \cite[Section 3.3.4]{LeThi_2018_30_years_dev} and \cite[Proposition 3]{rotaru2025tightanalysisdifferenceofconvexalgorithm}. Consequently, all the rates developed for DCA readily translate to rates for PGD using the mapping $\mu_1 = \gamma^{-1} + \mu_h$, $L_1 = \gamma^{-1} + L_h$, $\mu_2 = \gamma^{-1} - L_{\varphi}$, $L_2 = \gamma^{-1} - \mu_{\varphi}$. The condition $\mu_1 + \mu_2 > 0$ translates to the standard upper limit of $\gamma < \frac{2}{L_{\varphi} - \mu_{h}}$ and $\mu_2 < 0$ corresponds to stepsizes $\gamma > \frac{1}{L_{\varphi}}$.

When $\varphi$ is convex, \cite{Taylor_Jota_PGM_rates_proofs} shows various tight rates employing different performance metrics, including gradient residual for stepsize $\gamma = \frac{1}{L_{\varphi}}$, but not the gradient mapping. In the gradient residual, for $\varphi$ nonconvex, with $\mu_{\varphi} = -L_{\varphi}$ and $\mu_h=0$, $L_h=\infty$, the tight rate for stepsizes shorter than $\frac{1}{L_{\varphi}}$ is given in \cite[Proposition 8.11]{phd_abbaszadehpeivasti2024performance}. Sublinear rates for the same criterion are given in \cite{rotaru2025tightanalysisdifferenceofconvexalgorithm}. %

The projected gradient descent (ProjGD) is characterized by setting $h = \delta_{C}$ (indicator function of a non-empty, closed and convex set $C$), hence $\mu_h = 0$ and $L_h = \infty$. The proximal point algorithm is obtained by setting $\varphi$ constant.
The gradient descent (GD) is obtained by setting $h = 0$, hence $\mu_h = L_h = 0$. 
In this case, we further have $x^k - x^{k+1} = -\gamma \nabla \varphi(x^k)$, where $\gamma^{-1} = \mu_1$. Convergence rates on the best gradient norm are proved in \cite{rotaru2024exactworstcaseconvergencerates} for all curvatures of $\varphi$. These rates are a particular case of \cref{conjecture:Regimes_p5_above_thr} and served as a base to conjecture the tight expressions, since the projection / proximal step is not modifying the behavior of the worst-case function. 

To get exact rates from DCA, one just has to plug in the expressions of the curvatures $\mu_i$, $L_i$, $i=\{1,2\}$, written with respect to curvatures of $\varphi$ and $h$, in \cref{thm:dca_rates_N_steps,conjecture:Regimes_p5_above_thr}. The convexity of $h$ implies $\mu_h \geq 0$, thus $\mu_1 > 0$. Then regime $p_1$ corresponds to taking stepsizes $\gamma \leq \frac{1}{L_{\varphi}}$ for a splitting with $\varphi$ and $F$ nonconvex. Within the regimes $p_3$ and $p_4$, the stepsize $\gamma$ lies in the interval $[\frac{1}{L_{\varphi}}, \bar{\gamma}^{1}(L_{\varphi},\mu_{\varphi},\mu_h)]$, where $\bar{\gamma}^{1}(L_{\varphi},\mu_{\varphi},\mu_h)$ is a threshold with an analytical expression obtained through the condition $\mu_2 = \frac{-L_2 \mu_1}{L_2 + \mu_1}$. %
Regime $p_3$ holds when $\varphi$ and $F$ are nonconvex. When $\mu_h = 0$ and $\mu_{\varphi} = -L_{\varphi}$, then $\bar{\gamma}^{1}=\frac{\sqrt{3}}{L_{\varphi}}$, as also obtained in \cite{abbaszadehpeivasti2021GM_smooth} for the GD method. %
Regime $p_4$ holds when $F$ is strongly convex and $\mu_{\varphi} + \mu_{h} \geq 0$; in particular, this includes the case of $\varphi$ being strongly convex. For stepsizes $\gamma \in (\bar{\gamma}^{1}(L_{\varphi},\mu_{\varphi},\mu_h), \frac{2}{L_{\varphi}+\mu_{\varphi}})$, regime $p_5$ is active if $F$ is nonconcave. %
Regimes $p_2$ and $p_6$, respectively, correspond to using stepsizes $\gamma \leq \frac{1}{L_{\varphi}}$ for $\varphi$ concave and $F$ nonconcave or concave, respectively.

\cref{lemma:descent_lemma_stepsizes_geq_1} provides an example of convergence rate translation between PGD and DCA, corresponding to regime $p_3$ from \cref{lemma:regime_p3}, assuming the standard setup of PGD with $\mu_h = 0$ and $L_h = \infty$.%
\begin{lemma}\label{lemma:descent_lemma_stepsizes_geq_1}
    Let $\varphi \in \mathcal{F}_{\mu_{\varphi},L_{\varphi}}$, $h \in \mathcal{F}_{0,\infty}$ and $\gamma \in [\frac{1}{L_{\varphi}}, \frac{2}{L_{\varphi}})$. Then after two iterations of PGD with stepsize $\gamma$ connecting $x^{k}$, $x^{k+1}$ and $x^{k+2}$ it holds \[
     \Delta F(x^k)
        {}\geq{}
    \tfrac{[(2-L_{\varphi})(2-\mu_{\varphi}) - 1]}{\gamma(2-L_{\varphi}-\mu_{\varphi})} \tfrac{1}{2} \|\Delta x^{k}\|^2  {}+{} 
    \tfrac{1}{\gamma(2-L_{\varphi}-\mu_{\varphi})} \tfrac{1}{2} \|\Delta x^{k+1}\|^2.
\]
\end{lemma}
\cref{lemma:descent_lemma_stepsizes_geq_1} is the extension to PGD of the same result obtained for GD in \cite[Lemma 4.3]{rotaru2025tightanalysisdifferenceofconvexalgorithm}. Therein, it is involved in deriving a stepsize schedule gradually increasing from $\bar{\gamma}^{1}(L_{\varphi},\mu_{\varphi},0)$ towards $\frac{2}{L_\varphi + \mu_{\varphi}}$ (if $\varphi$ is convex), which gives a better worst-case guarantee than the best constant stepsize. Based on \cref{lemma:descent_lemma_stepsizes_geq_1}, the same schedule should give a similar performance for PGD.  Further, we show how it may improve the DCA convergence.%
%
\subsection*{Curvature shifting and schedules for DCA}
\begin{proposition}\label{prop:curvature_shifting_PGD}
    One PGD iteration with stepsize $\gamma$ on the splitting $F=\varphi + h$ produces the same iterate as one iteration with stepsize $\tilde{\gamma}$ on $F = \tilde{\varphi} + \tilde{h}$, with $\tilde{\varphi} = \varphi + \lambda \frac{\|\cdot\|^2}{2}$, $\tilde{h} = h - \lambda \frac{\|\cdot\|^2}{2}$, $\tilde{\gamma}^{-1} = \gamma^{-1} - \lambda$, where $\lambda < \gamma^{-1}$.
\end{proposition}
This relates to curvature shifting in DCA: $f^{\lambda}_i \coloneqq f_i - \frac{\lambda \|\cdot\|^2}{2}$, $i=\{1,2\}$, where $\lambda$ is a scalar. Based on a single-iteration analysis, \cite{rotaru2025tightanalysisdifferenceofconvexalgorithm} studies  the best curvature shifting, observing that it usually reduces the convexity of both functions. %
The procedure can be improved, as (optimized) stepsize schedules $\{\gamma^{k}(L_{\varphi},\mu_{\varphi})\}$ for PGD can be translated to DCA. Consider the standard setting of PGD with $L_\varphi > 0$, $\mu_h=0$, $L_h=\infty$, assume $\mu_1 > 0$ and choose $\gamma = \mu_1^{-1}$. 
Using \cref{prop:curvature_shifting_PGD}, a curvature splitting schedule for DCA is defined as $\{\lambda^k\} = \mu_1 - \{\gamma^{k}(\mu_1 - L_2, \mu_1 - \mu_2)\}$. 
%
%
\section{Conclusion}
We studied rates in the gradient mapping for DCA and PGD. Examining the behavior of DCA after two iterations, we identified six regimes partitioning the parameter space, and  proved the corresponding six tight sublinear convergence rates. Three of those sublinear rates hold for any number of iterations and, for the other regimes, we conjectured their exact convergence rates. We also show how PGD rates follow directly from DCA ones. This connection suggests the potential for an optimized splitting of the objective function used in DCA, possibly varying at each iteration, inspired by stepsize schedules used for (proximal) gradient descent. We leave the practical testing of those as future work.%
\bibliographystyle{plain}
\bibliography{ms.bib} 
\end{document}